\documentclass{article}
\usepackage{amsthm,amsmath,amssymb}

\newcommand{\nc}{\newcommand}

\newtheorem{thm}{Theorem}[section]

\newtheorem{lemma}[thm]{Lemma}

\newtheorem{corollary}[thm]{Corollary}

\newtheorem{definition}[thm]{Definition}

\nc{\Ext}{\operatorname{Ext}}
\nc{\FS}{\operatorname{FS}}
\nc{\NS}{\operatorname{NS}}
\nc{\Amp}{\operatorname{Amp}}
\nc{\Pic}{\operatorname{Pic}}
\nc{\Kom}{\operatorname{Kom}}
\nc{\DGrB}{\operatorname{DGrB}}
\nc{\antidiag}{\operatorname{antidiag}}
\nc{\diag}{\operatorname{diag}}
\nc{\mo}{\operatorname{mod}}
\nc{\Gr}{\operatorname{Gr}}
\nc{\Rep}{\operatorname{Rep}}
\nc{\Perf}{\operatorname{Perf}}
\nc{\Hom}{\operatorname{Hom}}
\nc{\Sym}{\operatorname{Sym}}
\nc{\RHom}{R\operatorname{Hom}}
\nc{\cRHom}{\operatorname{\mathcal{R}\mathcal{H}om}}
\nc{\cHom}{\operatorname{\mathcal{H}om}}
\nc{\End}{\operatorname{End}}
\nc{\Coh}{\operatorname{Coh}}
\nc{\Aut}{\operatorname{Aut}}
\nc{\Td}{\operatorname{Td}}
\nc{\Coker}{\operatornamoe{Coker}}
\nc{\coker}{\operatorname{coker}}
\nc{\colim}{\operatorname{colim}}
\nc{\Ker}{\operatorname{Ker}}
\nc{\img}{\operatorname{Im}}
\nc{\D}{\operatorname{D}}
\nc{\ch}{\operatorname{ch}}
\nc{\Stab}{\operatorname{Stab}}
\nc{\SL}{\operatorname{SL}}
\nc{\rk}{\operatorname{rk}}
\nc{\GL}{\operatorname{GL}}
\nc{\Log}{\mathop{\mathrm{Log}}}
\nc{\abs}[1]{\lvert#1\rvert}
\nc{\Cone}{\operatorname{Cone}}
\nc{\id}{\operatorname{id}}
\nc{\Li}{\operatorname{Li}}

\newcommand{\Db}{{\mathrm D}^{b}}

\nc{\cA}{{\mathcal A}}
\nc{\cB}{{\mathcal B}}
\nc{\cC}{{\mathcal C}}
\nc{\cD}{{\mathcal D}}
\nc{\cE}{{\mathcal E}}
\nc{\cF}{{\mathcal F}}
\nc{\cG}{{\mathcal G}}
\nc{\cH}{{\mathcal H}}
\nc{\cI}{{\mathcal I}}
\nc{\cJ}{{\mathcal J}}
\nc{\cK}{{\mathcal K}}
\nc{\cL}{{\mathcal L}}
\nc{\cM}{{\mathcal M}}
\nc{\cN}{{\mathcal N}}
\nc{\cO}{{\mathcal O}}
\nc{\cP}{{\mathcal P}}
\nc{\cQ}{{\mathcal Q}}
\nc{\cR}{{\mathcal R}}
\nc{\cS}{{\mathcal S}}
\nc{\cT}{{\mathcal T}}
\nc{\cU}{{\mathcal U}}
\nc{\cV}{{\mathcal V}}
\nc{\cW}{{\mathcal W}}
\nc{\cX}{{\mathcal X}}
\nc{\cY}{{\mathcal Y}}
\nc{\cZ}{{\mathcal Z}}

\nc{\bA}{{\mathbb A}}
\nc{\bB}{{\mathbb B}}
\nc{\bC}{{\mathbb C}}
\nc{\bD}{{\mathbb D}}
\nc{\bE}{{\mathbb E}}
\nc{\bF}{{\mathbb F}}
\nc{\bG}{{\mathbb G}}
\nc{\bH}{{\mathbb H}}
\nc{\bI}{{\mathbb I}}
\nc{\bJ}{{\mathbb J}}
\nc{\bK}{{\mathbb K}}
\nc{\bL}{{\mathbb L}}
\nc{\bM}{{\mathbb M}}
\nc{\bN}{{\mathbb N}}
\nc{\bO}{{\mathbb O}}
\nc{\bP}{{\mathbb P}}
\nc{\bQ}{{\mathbb Q}}
\nc{\bR}{{\mathbb R}}
\nc{\bS}{{\mathbb S}}
\nc{\bT}{{\mathbb T}}
\nc{\bU}{{\mathbb U}}
\nc{\bV}{{\mathbb V}}
\nc{\bW}{{\mathbb W}}
\nc{\bX}{{\mathbb X}}
\nc{\bY}{{\mathbb Y}}
\nc{\bZ}{{\mathbb Z}}

\begin{document}
\title{Homological mirror symmetry of Fermat polynomials}
\author{So Okada\footnote{Supported by JSPS Grant-in-Aid
\#21840030 and Kyoto University Global Center Of Excellence
Program; Email: okada@kurims.kyoto-u.ac.jp; Address:
Research Institute for Mathematical Sciences, Kyoto
University, 606-8502 Kyoto Japan.}}  \maketitle

\begin{abstract}
 We discuss homological mirror symmetry of Fermat polynomials in terms
 of derived Morita equivalence between derived categories of coherent
 sheaves and Fukaya-Seidel categories (a.k.a.  perfect derived
 categories \cite{Kon09} of directed Fukaya categories
 \cite{Sei08,Sei03,AurKatOrl08}), and some related aspects such as
 stability conditions, (kinds of) modular forms, and Hochschild
 homologies.
\end{abstract}

 \section{Introduction}
 Homological mirror symmetry was introduced by Kontsevich \cite{Kon95}
 as the categorical equivalence of so-called A- and B- models of
 topological field theories.  For a pair of manifolds $X$ and $Y$
 possibly with group actions, homological mirror symmetry is an
 equivalence of triangulated categories of coherent sheaves of $X$ and
 Lagrangians of $Y$ with natural enhancements to differential graded
 (dg) categories (see \cite{Kel06}).
  
  We have seen such equivalences for elliptic curves \cite{PolZas}, the
  quartic surface case \cite{Sei03}, degenerating families of Calabi-Yau
  varieties and abelian varieties \cite{KonSoi01}.  The framework has
  been extended to Fano varieties and singularities in
  \cite{AurKatOrl08, AurKatOrl06,Sei01}.  See \cite{HKKPTVVZ} for a
  comprehensive source of references.

  Let $X_{n}$ be the function $F_{n}:\bC^{n}\to \bC$ for the Fermat
  polynomial $F_{n}:=x_{1}^{n}+\cdots+x_{n}^{n}$, and $Y_{n}$ be
  $F_{n}:\bC^{n}/G_{n}\to\bC$ with the group
  $G_{n}:=\{(\xi_{k})_{k=1\ldots n}\in \bC^{n} \mid \xi_{k}^{n}=1\}$
  acting on $\bC^{n}$ as $(x_{k})_{k=1\ldots n}\in \bC^{n}\mapsto
  (\xi_{k}x_{k})_{k=1\ldots n}\in \bC^{n}$.  For Morsifications of our
  functions, we define Fukaya-Seidel categories as perfect derived
  categories of $A_{\infty}$ categories (see \cite{Kel01}) of ordered
  vanishing cycles (Lagrangian spheres) \cite{Sei08,Sei03,AurKatOrl08}.
  
  For $X_{n}$, we define $\Db(\Coh X_{n})$ as the perfect derived
  category of coherent sheaves on its fiber over zero in the projective
  space of $\bC^{n}$ (so, it gives the bounded derived category of the
  coherent sheaves of the Fermat hypersurface in $\bP^{n-1}$).  By the
  same logic, we define $\Db(\Coh Y_{n})$, which gives the bounded
  derived category of coherent sheaves on the Fermat hypersurface in
  $\bP^{n-1}$ as the stack with respect to the action of
  $H_{n}:=G_{n}/\mbox{diagonals}$.

  Between derived categories such as $\Db(\Coh Y_{n})$ and $\FS(X_{n})$,
  we prove derived Morita equivalences in the sense that we take compact
  generating objects on both sides and dg categories (derived
  endomorphisms) of the objects in some dg enhancements, and find them
  isomorphic.
  
  Homological mirror symmetry comes with the canonical mirror
  equivalence.  We put $\FS(Y_{n})$, as a quotient of $\FS(X_{n})$ by
  the dual group $\hat{H}_{n}$ of $H_{n}$. In $\FS(X_{n})$, the
  $A_{\infty}$ category of ordered vanishing cycles for a Morsification
  of $X_{n}$ is $A_{\infty}$ isomorphic to a dg category. We take the dg
  category of a compact generating object consisting of these vanishing
  cycles such that $\hat{H}_{n}$ acts as an automorphism group on the
  category of the dg category.  We put Fukaya-Seidel category
  $\FS(Y_{n})$ as the perfect derived category of the dg
  $\hat{H}_{n}$-orbit category \cite{Kel05} of the dg category.
  
  We take the dg category of a $H_{n}$-invariant compact generating
  object of $\Db(\Coh X_n)$ in some dg enhancement  and find the derived Morita equivalence between
  $\Db(\Coh X_{n})$ and $\FS(Y_{n})$, comparing the dg
  $\hat{H}_{n}$-orbit category and the dg category. The following is a
  summary:
  \begin{align*}
   \begin{array}{ccc}
    \Db(\Coh Y_{n}) & \cong &
     \FS(X_{n})\\
    \uparrow\mbox{equivariance by $H_{n}$}  & & \downarrow\mbox{
     quotient by $\hat{H}_{n}$}\\
   \Db(\Coh X_{n}) &\cong &   \FS(Y_{n}).
   \end{array}
  \end{align*}
  For recent discussion on some aspects of homological mirror symmetry,
  see \cite{KapKreSch, Kon09} for references.

  \section{Derived Morita equivalence}\label{sec:Mor}
  Let us say equivalence instead of derived Morita equivalence.  Let
  $A_{n-1}$ be the Dynkin quiver of $A$ type with $n-1$ vertices and
  one-way arrows.  The perfect derived category $\Db(\mo A_{n-1})$ is
  equivalent to the category of graded B-branes $\DGrB(x^{n})$ (see
  \cite{Orl,KajSaiTak}), which has the natural dg enhancement in terms
  of graded matrix factorizations.

  Hochschild cohomology of the path algebra of $A_{n-1}$ is trivial
  \cite{Hap}, so is that of their tensor products.  For Auslander-Reiten
  transformation $\tau$ of $\DGrB(x^n)$ and the projective-simple module
  $E$ of $\mo A_{n-1}$, $\DGrB(x^n)$ is equivalent to the perfect
  derived category of the dg category
  $\cA_{n-1}:=\End^{*}(\oplus_{\mu=0\ldots n-1}\tau^{-\mu}(E))$, which
  is generated by degree-zero idempotents and, if any, degree-one
  nilpotents.

  For the function $x^{n}:\bC \to \bC$, the perfect derived category of
  $\cA_{n-1}$ is equivalent to $\FS(x^{n}):=\FS(x^{n}:\bC\to\bC)$
  \cite[2B]{Sei01}. Ordered objects $E, \ldots, \tau^{2-n}(E)$ represent
  ordered vanishing cycles for a Morsification of the function. Let us
  mention that they may represent ones for the Morsification of
  $x^{n}+y^2:\bC^{2}\to\bC $ in \cite{ACa} (see \cite[Section
  3]{Sei01}). The object $\tau^{1-n}(E)$ represents a connected sum of
  vanishing cycles.

  Up to $A_{\infty}$ isomorphisms, the tensor product of the graded
  commutative algebra $\cA_{n-1}^{\otimes n}$ admits no non-trivial
  $A_{\infty}$ structures.  K\"{u}nneth formula for vanishing cycles of
  Fukaya-Seidel categories on morphisms and compositions
  \cite[Proposition 6.3, Lemma 6.4]{AurKatOrl08} implies that the
  perfect derived category of $\cA_{n-1}^{\otimes n}$ is equivalent to a
  full subcategory of $\FS(X_n)$.  By the classical theory of
  singularity, vanishing cycles represented as objects of
  $\cA_{n-1}^{\otimes n}$ compose of a generating object of $\FS(X_{n})$
  (see \cite[Conjecture 1.3]{AurKatOrl08} and \cite{Oka,SebTho} for
  related studies).

  Let us explain an equivalence between $\Db(\Coh Y_{n})$ and $\Db(\mo
  A_{n-1})^{\otimes n}$.  By \cite{Orl}, $\Db(\Coh X_n)$ is equivalent
  to the category of graded B-branes $\DGrB(F_n)$, which has the natural
  dg enhancement in terms of graded matrix factorizations and the shift
  functor $\tilde{\tau}$ such that $\tilde{\tau}^{n}\cong [2]$.  Let
  $\Omega_{n-1}:=\oplus_{\mu=0\ldots
  n-1}\tilde{\tau}^{-\mu}(\cO_{X_{n}})$.  With explicit descriptions
  given by graded matrix factorizations for each summand of
  $\Omega_{n-1}$ (see \cite[Section IV]{Asp07}), we
  have the following. The object $\Omega_{n-1}$ is isomorphic to the
  compact generating object $\oplus_{\mu=0\ldots
  n-1}\Omega_{\bP^{n-1}}^{n-1-\mu}(n-\mu)[-\mu]|_{X_{n}}$ of $\Db(\Coh
  X_{n})$.  The object $\Omega_{n-1}$ is invariant under $H_{n}$ actions
  up to canonical isomorphisms $\hat{H}_{n}$, and so $H_{n}$ equivariant
  objects of $\Omega_{n-1}$ compose of a generating object of $\Db(\Coh
  Y_{n})$.  For the dg category of $\Omega_{n-1}$ in the natural dg
  enhancement of $\DGrB(X_{n})$, the $H_{n}$-equivariant dg category of
  the dg category of $\Omega_{n-1}$ coincides with $\cA_{n-1}^{\otimes
  n}$.
  
  Let us explain an equivalence between $\Db(\Coh X_{n})$ and
  $\FS(Y_{n})$.  Morphism spaces of the dg category of $\Omega_{n-1}$
  are representations of $H_{n}$ and compositions of morphisms are
  multiplications of matrices.  In the sense of \cite{CibMar}, by the
  action of $H_{n}$, morphism spaces of the dg category of
  $\Omega_{n-1}$ are $\hat{H}_{n}$-graded.  For a finite group $G$ of dg
  functors of a dg category $\cA$ inducing an equivalence on
  $H^{0}(\cA)$, we have the dg orbit category \cite{Kel05} $\cB=\cA/G$
  such that dg categories $\cA$ and $\cB$ have the same objects and for
  objects $X$ and $Y$ in $\cB$, $\cB (X,Y)=\colim_{g\in G}\oplus_{f\in
  G} \cA (f(X), g\circ f(Y))$.  By applying \cite[Proposition
  3.2]{CibMar}, the $\hat{H}_{n}$-orbit category of the category of
  $\cA_{n-1}^{\otimes n}$ is equivalent to the category of the dg
  category of $\Omega_{n-1}$.  Differentials have not been altered by
  taking dg equivariance and orbits. The dg category of $\Omega_{n-1}$
  in $\DGrB(X_{n})$ is equivalent to the dg $\hat{H}_{n}$-orbit category
  of $\cA_{n-1}^{\otimes n}$.  Explicitly, we define $\FS(Y_{n})$ as the
  perfect derived category of the dg $\hat{H}_{n}$-orbit category.  The
  perfect derived category of the dg category of $\Omega_{n-1}$ in
  $\DGrB(X_{n})$ is equivalent to $\Db(\Coh X_{n})$, and by the
  construction, $\FS(Y_{n})$ is equivalent to $\Db(\Coh X_{n})$. So, we
  have the following.

  \begin{thm}\label{thm:main}
   For $n>0$, $\Db(\Coh Y_n)$ is derived Morita equivalent to $\FS(X_n)$
   and $\Db(\Coh X_n)$ is derived Morita equivalent to $\FS(Y_n)$.
  \end{thm}

  \section{Discussion}

  Let us start off with some aspects of the dg category of
  $\Omega_{n-1}$, which is generated by degree-zero idempotents
  representing summands and, if any, degree-one closed differential
  forms.  The dg category of $\Omega_{n-1}$ together with its
  deformation theory would be related to the notion of Gepner model in
  physics (see \cite{Asp05,HerHorPag}). In \cite[Section 5.1]{HorWal},
  explicit forms of deformations of summands of the dg category of
  $\Omega_{4}$ have been computed along the deformation flow which is
  compatible with the actions of $(\overline{\xi}_{1},\ldots,
  \overline{\xi}_{5})$ in $H_{5}$ subject to $\overline{\xi}_{1}\cdots
  \overline{\xi}_{5}=1$.

  Each summand of $\Omega_{n-1}$ by the definition represents a
  Lagrangian among $\hat{H}_{n}$ isomorphic ones. In $\FS(x^{n})$, in
  the notation of the previous section, we observe that we count the
  intersection of the connected sum $\tau^{1-n}(E)$ with vanishing
  cycles $E$ and $\tau^{2-n}(E)$. We have non-trivial homomorphisms from
  $\tau^{1-n}(E)[1]$ to $E$ and from $\tau^{2-n}(E)$ to
  $\tau^{1-n}(E)[1]$ and trivial ones with $\tau^{-\mu}(E)$ for
  $0<\mu<n-2$.  For $0 \leq \mu<\mu'\leq n-1$ and Euler paring $\chi$,
  $\chi(\tilde{\tau}^{-\mu}(\cO_{X_{n}})[\mu],
  \tilde{\tau}^{-\mu'}(\cO_{X_{n}})[\mu'])$ counts intersections for
  each Lagrangian represented by $\tilde{\tau}^{-\mu}(\cO_{X_{n}})$ with
  Lagrangians represented by $\tilde{\tau}^{-\mu'}(\cO_{X_{n}})$.

  \subsection{Two types of stability conditions}
  Let us discuss the notion of stability conditions on triangulated
  categories \cite{BRI,GorKulRud,KonSoi08}. From a view point, a
  stability condition on a triangulated category is a collection of
  objects, which are called stable and has a partial order with indices
  called slopes, such that for consecutive stable objects $E \geq F$
  with a non-trivial homomorphism from $E$ to $F[1]$ we have consecutive
  stable objects $F[1]\geq E$.  For reasonable such collections, we have
  taken extension-closed full subcategories consisting of some
  consecutive stable objects to study wall-crossings, moduli problems,
  Donaldson-Thomas type invariants, and so on.

  For the quintic, a few types of stability conditions have been studied
  in mathematics and physics.  First type of stability conditions
  consists of variants of Gieseker-Maruyama stabilities such as
  \cite{Bay,Ina,Tod}.  Slopes of stable objects are given by appropriate
  approximations of Hilbert polynomials. These stability conditions are
  compatible with the tensoring by the ample line bundle, which is
  called the monodromy around the large radius limit (see \cite[Section
  F]{Asp07}).
  
  Second type of stability conditions has been studied in terms of
  Douglas' $\Pi$-stabilities for Gepner point \cite{AspDou,Dou02,
  Dou01}.  With our view point, the simplest stability condition
  compatible with $\tilde{\tau}$, which is called the monodromy around
  Gepner point (see \cite[Section F]{Asp07}), and with complexfied
  numerical classes of vanishing cycles, is as follows.  For $\mu\in
  \bZ$, stable objects consist of objects
  $\tilde{\tau}^{-\mu}(\cO_{X_{5}})$ and the order of them is given by
  the clock-wise order of $\exp(-\mu\frac{2\pi i}{5})$.

  The extension-closed full subcategory of the consecutive stable
  objects $\cO_{X_{5}}>\tilde{\tau}^{-1}(\cO_{X_{5}})$ of the stability
  condition of the second type gives the graded generalized Kronecker
  quiver with degree-one arrows from $\cO_{X_{5}}$ to the other,
  together with the same number of degree-two arrows in the opposite way
  and a degree-three loop on each vertex.  For some vanishing cycle
  represented by $\cO_{X_{5}}$, the number of degree-one arrows
  represents the number of intersections of the vanishing cycle with
  vanishing cycles represented by the other object and coincides with
  the Donaldson-Thomas invariant of the moduli space of the numerical
  class which is the sum of the ones of the two stable objects.  The one
  of the three consecutive stable objects
  $\cO_{X_{5}}>\tilde{\tau}^{-1}(\cO_{X_{5}})>\tilde{\tau}^{-2}(\cO_{X_{5}})$
  gives two copies of the quiver glued, together with degree-two arrows
  from $\cO_{X_{5}}$ to $\tilde{\tau}^{-2}(\cO_{X_{5}})$ and the same
  number of degree-one arrows in the opposite way and a degree-three
  loop on each vertex.  The objects $\cO_{X_{5}}$,
  $\tilde{\tau}^{-1}(\cO_{X_{5}})$, and $\tilde{\tau}^{-2}(\cO_{X_{5}})$
  make cluster collections \cite{KonSoi09, KonSoi08}.

  The stability condition of the second type is equipped with the
  consecutive stable objects
  $\tilde{\tau}^{-1}(\cO_{X_{5}})[1]>\cO_{X_{5}}$.  The extension-closed
  full subcategory of the pair gives the graded generalized Kronecker
  quiver with degree-zero arrows in one way, together with the same
  number of degree-three arrows in the opposite way and a degree-three
  loop on each vertex, which coincides with that of the consecutive
  stable objects $\cO_{X_{5}}(1)>\cO_{X_{5}}$ of the simplest stability
  condition of the first type on the sequence $\cO_{X_{5}}(\mu)$ for
  $\mu \in \bZ$.

  \subsection{(kinds of) Modular forms}
  For a stability condition of either type, objects in the
  extension-closed full subcategory consisting of stable objects of a
  slope are said to be semistable.  Let us count semistable objects in
  the extension-closed full subcategory consisting of a stable object for
  the stability condition of the second type.  The same counting logic
  applies when we have a single stable spherical object of some slope.

  As in \cite{MelOka}, we are interested in (kinds of) modular forms and
  Donaldson-Thomas type invariants.  For the motive of the affine line
  $\bL$, we represent the formal symbol $\bL^{\frac{1}{2}}$ by
  $-q^{\frac{1}{2}}$, which conjecturally corresponds to so-called
  graviphoton background in the supergravity theory
  \cite{DimGuk,DimGukSoi}.  We obtain $\cJ_{m}(q):=
  \frac{1}{m}\frac{q^{\frac{m}{2}}}{(1-q^{m})}$ for $m>0$ as the
  coefficient of $x^m$ in the formal logarithm of the quantum
  dilogarithm $\sum_{m\geq
  0}\frac{(-q^{\frac{1}{2}})^{m^2}}{(q^{m}-1)\cdots (q^{m}-q^{m-1})}x^m$
  \cite{KonSoi09,KonSoi08}.  The generating function
  $\sum_{n>0}m^{1-k}\cJ_{m}(q)= \sum_{m,r>0}
  m^{-k}q^{\frac{mr}{2}}-m^{-k}q^{mr}$ consists of Eisenstein series for
  negative odd $k$, or Eichler integrals \cite{LawZag} of Eisenstein
  series for positive odd $k$.  In particular, for $k=-1$, $q=\exp(2\pi i
  z)$, and Eisenstein series $E_{2}(z)$, we obtain the quasimodular form
  $E_{2}(z)-E_{2}(2z)$ \cite{KanZag}.

  \subsection{Hochschild homologies}
  In $\bP^{n-1}$, $Y_{n}$ is the same as $\bP^{n-2}$ given by
  $x_{1}+\cdots+x_{n}=0$ with the orbifold structures of degree $n$
  along coordinate hypersurfaces.  By computing Poincar\'{e} polynomial
  of $Y_{n}$, which can be partitioned according to degrees of orbifold
  structures, we obtain $P(Y_{n},q)= \sum_{2 \leq j \leq n, 2 \leq k
  \leq j} n^{n-j} \binom{n}{j} \ (-1)^{j-k} P(\bP^{k-2},q)$.
  For example, we have $1$,
  $q+7$, $q^2+13 q+67$, and $q^3+21 q^2+181q+821$.  In particular, by
  computing Euler characteristics in terms of Hochschild homologies (see
  \cite{Kel98, HocKosRos}), we obtain $P(Y_n,1)=(n-1)^{n}$.
 
 \section{Acknowledgments}
 The author thanks Australian National University, Institut des Hautes
 \'Etudes Scientifiques, and Research Institute for Mathematical
 Sciences for providing him postdoctoral support.  He thanks Professors
 Bridgeland, Fukaya, Hori, Katzarkov, Keller, Kimura, Kontsevich,
 Lazaroiu, Mizuno, Moriyama, Nakajima, Neeman, Polishchuk, Seki,
 A. Takahashi, Terasoma, Usnich, Yagi, and Zagier for their useful
 discussions. In particular, he thanks Professor Kontsevich for his
 generous and inspiring discussions.

\end{document}